\documentclass{amsart}

\usepackage{amsmath, amsthm, amssymb}
\usepackage{epsfig}
\usepackage{graphicx}
\usepackage{array}
\usepackage{amsmath}
\usepackage{amsfonts}
\usepackage[all]{xy}
\usepackage{amsfonts}
\usepackage{amsmath}%
\usepackage{amsfonts}%
\usepackage{amssymb}%
\usepackage{graphicx}

\usepackage{amsfonts}
\usepackage{amssymb}
\usepackage{amsmath,amsxtra,amssymb}
\numberwithin{equation}{section}

\makeatletter
\mathchardef\hugecheck="7014
\newcommand\hugesize{\@setfontsize\hugesize{25pt}{0}}
\newcommand\smallhugesize{\@setfontsize\smallhugesize{20pt}{0}}

\def\specialchecksmall{\mbox{\hbox
to0pt{\raisebox{-4pt}{\smallhugesize$\hugecheck$}}
                                 $\kern-2.7pt\otimes$}}
\makeatother

\newtheorem{theorem}{Theorem}[section]
\newtheorem{lemma}[theorem]{Lemma}
\newtheorem{remark}[theorem]{Remark}
\newtheorem{corollary}[theorem]{Corollary}
\newtheorem{proposition}[theorem]{Proposition}

\newcommand{\field}[1]{\mathbb{#1}}
\newcommand{\C}{{\field{C}}}

\newcommand{\ra}{\rightarrow}
\newcommand{\id}{{\iota}}
\newcommand{\ot}{{\otimes}}
\newcommand{\om}{{\omega}}
\newcommand{\tp}{{\widehat{\otimes}}}
\newcommand{\vtp}{{\overline{\otimes}}}

\newcommand{\h}{{\mathcal H}}

\newcommand{\T}{{\mathcal{T}}}
\newcommand{\B}{{\mathcal{B}}}

\newcommand{\fee}{{\varphi}}

\newcommand{\LL}{{L_{\infty}(\G)}}
\newcommand{\LO}{{L_{1}(\G)}}
\newcommand{\LT}{{L_{2}(\G)}}
\newcommand{\MG}{{M(\G)}}
\newcommand{\CZ}{{C_0(\G)}}
\newcommand{\CU}{{C_u(\G)}}

\newcommand{\loneqg}{L_1(\mathbb{G})}
\newcommand{\linfqg}{L_\infty(\mathbb{G})}

\newcommand{\PU}{\mathcal {P}}
\newcommand{\G}{\mathbb G}

\def\proclaim #1. #2\par{\medbreak
\noindent{\bf#1.\enspace}{\sl#2}\par\medbreak}
\allowdisplaybreaks

\begin{document}

\title[Poisson Boundaries]
{Poisson Boundaries over Locally Compact Quantum  Groups}
\author{Mehrdad Kalantar, Matthias Neufang and Zhong-Jin Ruan}
\address{ School of Mathematics and Statistics,
         Carleton University, Ottawa, Ontario, Canada K1S 5B6}
\email{mkalanta@math.carleton.ca}
\address{School of Mathematics and Statistics,
         Carleton University, Ottawa, Ontario, Canada K1S 5B6}
\email{mneufang@math.carleton.ca}
\address{Universit\'{e} Lille 1 - Sciences et Technologies,
UFR de Math\'{e}matiques, Laboratoire de Math\'{e}matiques Paul Painlev\'{e} - UMR CNRS 8524,
59655 Villeneuve d'Ascq C\'{e}dex, France}
\email{Matthias.Neufang@math.univ-lille1.fr}


\address{Department of Mathematics,
         University of Illinois, Urbana, IL 61801, USA}
\email{ruan@math.uiuc.edu}
\thanks{The second and third author were partially supported by NSERC, 
and the National Science
Foundation DMS-0901395, respectively.}

\subjclass[2000]{Primary 46L53, 46L89; Secondary 46L07, 46L65, 60J50.}

\begin{abstract}
We present versions of several classical results on harmonic functions and Poisson boundaries in the setting of locally compact quantum groups. 
In particular, the Choquet--Deny theorem holds for compact quantum groups;
also, the result of Kaimanovich--Vershik and Rosenblatt,
which characterizes group amenability in terms of harmonic functions, 
admits a non-commutative analogue in the separable case.
We also explore the relation between classical and quantum Poisson boundaries
by investigating the spectrum of the quantum group.
We apply this machinery to find a concrete realization of the Poisson boundaries of 
the compact quantum group $SU_{q}(2)$  arising from measures on its spectrum.
\end{abstract}


\maketitle

\section{Introduction}

The theory of von Neumann algebras originated 
in a series of remarkable papers during the late 1930s and early 1940s
by Murray and von Neumann.
The theory may be viewed as an operator, or noncommutative, version of measure theory.
During the last seventy years, operator algebras have proved to have a 
very profound structure theory.
They also provide the foundation to consider the quantization of 
many areas of mathematics, such as analysis, topology, geometry,
probability, and ergodic theory.
Recently, work of Woronowicz, Baaj--Skandalis, 
and Kustermans--Vaes, has led to the very successful development of the theory of locally compact quantum groups.
This provides the natural framework for  the quantization of various problems related to groups and 
group actions on measure (or topological) spaces.
The aim of this paper is to study Poisson boundaries
over locally compact quantum groups.

Poisson boundaries and harmonic functions have played a very important role in the study of random
walks on discrete groups, and more generally in harmonic analysis and
ergodic theory on locally compact groups
(see for instance  Furstenberg's seminal work \cite{Furst}).
Let us recall that if $G$ is a locally compact group and  
$\mu$ is  a probability measure on $G$, we
obtain a Markov operator $\Phi_{\mu}$ on  $L_{\infty}(G)$ 
associated with the measure $\mu$ which is 
defined by
\begin{equation}
\label {F.1.1}
\Phi_{\mu}(h) (s)=  \mu \star h(s)= \int_{G} h(st){d}\mu(t) \ \ \ \ (s\in G).
\end{equation}
It is known that  there exists 
a probability measure space $(\Pi, \nu)$, the Poisson boundary of $(G, \mu)$,
such that $L_\infty(\Pi, \nu)$ can be identified with  the weak* closed subspace 
$\h_\mu$ of $L_{\infty}(G)$ which consists of all 
 $\mu$-harmonic functions, i.e., functions $h$ on $G$ satisfying $\Phi_{\mu}(h) = h$.

Poisson boundaries over discrete quantum groups
 were first studied by Izumi \cite{I}.
He used these objects to study compact quantum group actions. 
 More precisely, he showed that the relative
 commutant of the fixed point algebra of certain ITP actions 
 of compact quantum groups can be realized as the
 Poisson boundary over the dual (discrete) quantum group.
 As a concrete example, he showed that in the case of Woronowicz' compact quantum group 
 $SU_q(2)$, this gives an identification between the  Poisson boundary
 of a Markov operator on the dual quantum group, with one of
 the Podle\'{s}  spheres \cite{Podles}. 
Later this result was extended to the
case of $SU_q(n)$ by Izumi, Neshveyev and Tuset \cite{INT}, and further
generalized by Tomatsu in \cite{Tom}.
Poisson boundaries for other discrete quantum groups have been studied
by Vaes, Vander Vennet and Vergnioux  \cite{VV1}, \cite{VV2}, \cite{VV3}.

In this paper, we establish important classical results
on Poisson boundaries and (bounded) harmonic functions
in the general quantum group setting.
The paper is organized as follows.
We recall relevant definitions and introduce some notation in section 2.
In section 3, we establish quantized versions of several classical results
concerning Poisson boundaries.
In particular, we prove that the  Poisson boundary of a non-degenerate
`{\it quantum probability measure}' on a locally compact quantum group is never a
subalgebra, unless trivial, and that there is no non-trivial harmonic operator
which is `{\it continuous}' and `{\it vanishing at infinity}'.

\par
In the classical setting there is a characterization of amenability of a ($\sigma$-compact)
locally compact group in terms of its Poisson boundaries. 
Kaimanovich--Vershik \cite{1} and Rosenblatt \cite{Ros} independently proved that
if $G$ is a ($\sigma$-compact) locally compact amenable group, then there exists
an absolutely continuous measure $\mu$ on $G$ (i.e., $\mu\in L_{1}(G)$)
such that  $\mu$-harmonic functions are trivial. 
This answered a conjecture of Furstenberg \cite{Furst},
who had shown the converse.
 In section 4 we prove the corresponding result in the quantum setting (Theorem \ref{KV}).

\par
The classical Choquet--Deny theorem \cite{CD} states that there is no non-trivial $\mu$-harmonic function
 for an adapted probability measure on a locally compact abelian group $G$.
 The conclusion of this theorem has been proved for many other cases,
 including compact groups. In section 5 we study 
 Poisson boundaries over compact quantum groups, and we prove a noncommutative
 version of the Choquet--Deny theorem in this setting (Theorem \ref{cdthem}).

\par
We investigate the relation between the classical and the quantum setting in section 6,
 by proving a formula which links the  Poisson boundary of a Markov 
 operator induced
 from a commutative quantum subgroup, to its classical counterpart. Applying our machinery
 to the case of $SU_q(2)$, we show that the Poisson boundary over the latter,
 induced from the quantum subgroup $\mathbb{T}$, i.e., the spectrum
 of $SU_q(2)$, can be identified with a Podle\'{s}  sphere.


\section{Definitions and preliminary results}
In this paper, we denote by   $\G = (\LL, \Gamma , \varphi , \psi )$ 
 a (von Neumann algebraic) \emph{locally compact quantum group}
in the sense of  Kustermans and Vaes  \cite{KV1}, \cite{KV2}.
 The  \emph {right Haar weight} $\psi$ (which is an n.s.f. right invariant weight)
  on the quantum group von Neumann algebra 
 $\LL$ determines a Hilbert space  $L_{2}(\G)= L_{2}(\G, \psi)$ and 
 we obtain   the \emph{right fundamental unitary operator}  $V$
on $L_{2}(\G) \otimes L_{2}(\G)$, which satisfies the
  \emph{pentagonal relation} 
 \begin{equation}
\label {F.pentagonal} V_{12} V_{13}V_{23} = V_{23} V_{12}.
\end{equation}
Here we used the leg notation
$V_{12}=V\ot 1$, $V_{23} = 1\ot V$, and $V_{13} = (\id\ot\chi)V_{12}$,
where $\chi(x\ot y) = y\ot x$ is the flip map.
This fundamental unitary operator induces 
a  coassociative comultiplication 
\begin{equation}
\label {F.general}
\tilde \Gamma: x \in \B(L_{2}(\G)) \to \tilde \Gamma(x) = V(x\otimes 1)V^{*}
\in \B(L_{2}(\G)\otimes L_{2}(\G))
\end{equation}
on $\B(L_{2}(\G))$, for which   we have 
$\tilde \Gamma _{|L_{\infty}(\G)}= \Gamma$.

Let $\loneqg$ be the predual of $\linfqg$. Then the pre-adjoint of
$\Gamma$ induces  an associative completely
contractive multiplication
\begin{equation}
\label {F.mul} \star  :  f \otimes g \in 
L_{1}(\G)\tp L_{1}(\G) \to f \star g = (f \otimes
g) \Gamma \in L_{1}(\G)
\end{equation}
on $\loneqg$.
Since the multiplication $\star$ is a complete quotient map from
$L_{1}(\G)\hat \otimes L_{1}(\G)$ onto $L_{1}(\G)$, we get 
\begin{equation}
\label {F.L2} \loneqg = \langle  L_{1}(\G) \star L_{1}(\G) \rangle  =
\overline{\mbox{span}\{f\star g: f, g\in \loneqg\}}^{\|\cdot\|}.
\end{equation}
If $\G_{a}= (L_{\infty}(G), \Gamma_{a}, \varphi_{a}, \psi_{a})$ 
is the  commutative quantum group associated with a locally compact group $G$, 
then $L_{1}(\G_{a})$ is just the convolution algebra $L_{1}(G)$.
If on the other hand $\G_{s}=\hat \G_{a}$ is the cocommutative 
dual quantum group of $\G_{a}$, then
 $L_{1}(\G_{s})$  is the Fourier algebra $A(G)$.

The  \emph{right  regular representation} $\rho : L_{1}(\G) \to
\B(L_{2} (\G))$ is defined by
 \[
\rho : f \in L_{1}(\G)   \to \rho(f) = (\id\otimes f)(V)
\in \B(L_{2}(\G)),
 \]
which is an injective and completely contractive algebra homomorphism
from $L_{1}(\G)$ into $\B(L_{2} (\G))$.
We let  $L_{\infty} (\hat \G')={\{\rho(f): f\in \loneqg\}}''$
 denote  the quantum group von Neumann algebra 
of  the (commutant) dual quantum group $\hat \G'$.
Then   $\hat V = \Sigma V^{*}\Sigma$,
where $\Sigma$ denotes the \emph{flip operator} $\Sigma (\xi \otimes \eta)
= \eta \otimes \xi$ on $L_{2}(\G)\otimes L_{2}(\G)$,
is the  right fundamental unitary operator   of $\hat \G'$, and   
\begin{equation}
\label {F.hatrg}
\hat\rho: \hat f' \in L_{1}(\hat \G ') \to \hat\rho(\hat f' ) =
(\id \ot \hat f')(\hat V) = (\hat f ' \ot \id )(V^{*})\in \LL
\end{equation}
is the  {right  regular representation} of $\hat \G'$.
The \emph{reduced quantum group $C^*$-algebra}  
 $C_{0}(\G) = \overline{\hat\rho(L_{1}(\hat \G') )}^{\|\cdot\|}$
  is a weak$^*$ dense $C^*$-subalgebra of $\LL$ with the 
 comultiplication 
\[
\Gamma :  C_{0}(\G)\to M(C_{0}(\G)\otimes
C_{0}(\G))
\]
given by the  restriction of   the comultiplication of $\LL$  to $C_{0}(\G)$.
Here,  $M(C_{0}(\G)\otimes C_{0}(\G))$  denotes the multiplier $C^*$-algebra of
the minimal $C^*$-algebra tensor product $C_{0}(\G)\otimes
C_{0}(\G)$. For convenience, we often use $C(\G)$ for $M(C_{0}(\G))$.
Let $M(\G)$ denote the operator dual $C_{0}(\G)^{*}$.
There exists  a completely contractive multiplication on $M(\G)$ given by
the convolution
 \[
 \star : \mu\ot \nu \in M(\G)\tp M(\G)\mapsto \mu \star \nu 
 = \mu (\id\otimes \nu)\Gamma = \nu (\mu \otimes \id)\Gamma
 \in M(\G)
 \]
such that  $M(\G)$ contains $\loneqg$ as a norm closed two-sided ideal
(for details see \cite {BT} and \cite {HNR}).

Let $L_{1 *}(\hat\G') = \{\hat\om'\in L_{1}(\hat \G') : \exists \hat f'\in L_{1}(\hat \G') \ \text{such that}
 \ \hat\rho(\hat\om')^* = \hat\rho(\hat f')\}$.
Then $L_{1 *}(\hat\G')\subseteq L_{1}(\hat \G')$ is norm dense, and with the involution
$(\hat\om')^* = \hat f'$, and the norm $\|\hat\om'\|_u = max\{\|\hat\om'\| , \|(\hat\om')^*\|\}$, 
the space $L_{1 *}(\hat\G')$
becomes a Banach $^*$-algebra
(see \cite {Kus} for details).
We obtain the   \emph{universal quantum group $C^*$-algebra} $\CU$
as the universal enveloping $C^*$-algebra of the Banach algebra $L_{1 *}(\hat\G')$.
There is a universal $^*$-representation
\[
\hat \rho_{u} : L_{1 *}(\hat \G')\to \B(H_{u})
\]
such that $\CU = \overline {\hat \rho_{u}(L_{1}(\hat \G'))}^{\|\cdot\|}$.
 There is  a universal comultiplication
\[
\Gamma_u :\CU\ra M(\CU\otimes\CU), 
\]
and the operator dual $M_u(\G)=C_{u}(\G)^{*}$, which can be regarded
as the space of all {\it quantum measures} on $\G$, is a 
unital completely contractive Banach 
algebra with multiplication given by
\[
\om \star_{u} \mu = \om (\id\otimes \mu)\Gamma_{u} = \mu(\om\otimes \id)\Gamma_{u}
\]
(see \cite {BS}, \cite {BT} and \cite {Kus}).
By the universal property of $\CU$, there is a  unique surjective $*$-homomorphism 
$\pi:\CU\ra\CZ$  such that $\pi(\hat\rho_u(\hat\om')) = \hat\rho(\hat\om')$
for all $\hat\om'\in L_{1 *}(\hat \G')$. Moreover,
the adjoint map $\pi^{*}: M(\G)\to M_u(\G)$ defines a completely isometric injection
such that $\mu\star_{u}\pi^*(\om)$ and $\pi^{*}(\om)\star_{u} \mu$ are in $\pi^*(M(\G))$
for all $ \mu\in M_u(\G)$ and $\om\in\MG$  (\cite [Proposition 6.2]{Kus}).
Therefore we can identify $M(\G)$ with a
norm closed two-sided ideal in $M_u(\G)$, and
\begin{equation}
\label {F.2mul}
\mu\star \om = (\pi^*)^{-1}(\mu\star_{u}\pi^*(\om))\in\MG \ \ \ \ ,
\ \ \ \  \om \star \mu = (\pi^{*})^{-1}(\pi^{*}(\om)\star_{u} \mu)\in M(\G)
\end{equation}
define actions of $M_u(\G)$ on $M(\G)$.
In particular, the restriction of $\pi^{*}$ to $\LO$ is a completely isometric injection
from $\LO$ into $M_u(\G)$.  Since $\LO = \langle \LO\star \LO\rangle$,
we can conclude that $\mu \star f$ and $f \star \mu$ are again contained in
$\LO$.  Therefore, we can also identify $\LO$ with a norm closed two-sided ideal in $M_u(\G)$.
%
%
%
%
%
%
Then it is seen from (\ref {F.2mul}) that for each $\mu\in M_u(\G)$, we obtain
a pair of completely bounded maps
\begin{equation}
\label {F.multi}
\mathfrak{m}^{l}_{\mu}(f) = \mu \star f ~ \ \ \ \ \mbox{and } ~ \ \ \ \
\mathfrak{m}^{r}_{\mu}(f) = f \star \mu
\end{equation}
on $\LO$ with $\mbox{max}\{\|\mathfrak{m}^{l}_{\mu}\|_{cb},
\|\mathfrak{m}^{r}_{\mu}\|_{cb}\}\le \|\mu\|$.
%
%
%
The adjoint map $\Phi_{\mu}= (\mathfrak{m}^{r}_{\mu})^{*}$ 
is a normal completely bounded map  on $\LL$ such that $\Phi_{\mu}(x) = \mu\star x$,
more precisely, 
\begin{equation}
\label {F.210}
\langle f , \Phi_{\mu} (x)\rangle = \langle f \star \mu,  x \rangle 
= \langle f, \mu\star x\rangle
\end{equation}
for all $x \in \LL$ and $ f \in \LO$.
Moreover, the map $\Phi_\mu$ satisfies the covariance condition
\begin{equation}
\label {F.Phi}
\Gamma \circ \Phi_{\mu} = (\id \otimes \Phi_{\mu})\circ \Gamma, 
\end{equation}
or equivalently,
$\Phi_{\mu}(x \star f) = \Phi_{\mu}(x)\star f$
for all $x \in \LL$ and $f\in \LO$. 
If we let
\[
LUC(\G) = \langle \LL\star \LO\rangle = \overline {\mbox{span} \{
x\star f: x \in \LL, f\in \LO\}}^{\|\cdot\|}
\]
denote  the space of  \emph{left uniformly continuous linear functionals} on $\LO$, then 
$\Phi_{\mu}$  maps $LUC(\G)$ into $LUC(\G)$.
Under certain  conditions, $LUC(\G)$ is a  unital $C^*$-subalgebra of $C(\G)$ 
(cf. \cite {HNR2}, \cite {HNR} and \cite {Run2}). 
Also, since $\Phi_{\mu}$  maps $\CZ$ 
into $\CZ$, if $\Phi_{\mu}$ is completely positive, it maps $C(\G)$ into $C(\G)$
(see \cite {Lance}). 
In general, we have 
\[
C_{0}(\G)\subseteq LUC(\G)\subseteq C(\G)\subseteq \LL.
\]
In particular, if $\G$ is a compact quantum group,  we have 
$C_{0}(\G)=  LUC(\G)= C(\G)$,
and if $\G$ is discrete, we have $LUC(\G)= C(\G) = \LL$.



\section{Poisson boundaries of quantum probability measures}

We denote by $\PU_u(\G)$ the set of all states on $\CU$ (i.e., the `quantum probability measures'). 
Then $\Phi_\mu$ is a \emph{Markov operator}, i.e.,
a unital normal completely positive map, on $L_{\infty}(\G)$.

We consider the space of fixed points $\h^{\mu} = \{x\in L_{\infty}(\G): \Phi_{\mu}(x)=x\}$.
It is easy to see that $\h^{\mu}$ is a weak*
closed operator system in $L_{\infty}(\G)$.
In fact, we obtain a natural von Neumann algebra product on this space.
Let us  recall this construction
for the convenience of the reader (cf. \cite [Section 2.5]{I}).

We first  define a  projection
${\mathcal E}_{{\mu}} : L_{\infty}(\G) \to L_{\infty}(\G)$  of norm one by taking 
 the weak$^*$ limit
\begin{equation}\label{11}
{\mathcal E}_{\mu} (x) = \lim_{\mathcal{U}} \frac{1}{n}\sum_{k=1}^{n} 
\Phi_{\mu}^k (x) 
\end{equation}
with respect to   a   free ultrafilter  $\mathcal{U}$  on
$\mathbb{N}$.
It is easy to see that  $\h^{\mu} = {\mathcal E}_{\mu}(L_{\infty}(\G))$, and that
the Choi-Effros product 
\begin{equation}
\label {2l}
x\circ y  = {\mathcal E}_{\mu} (xy) 
\end{equation}
defines  a von Neumann algebra product on $\h^{{\mu}}$.
We note that this product is independent of the choice of the free ultrafilter 
$\mathcal{U}$ since every completely positive isometric linear isomorphism between
two von Neumann  algebras is a $^*$-isomorphism.
To avoid confusion, we denote by $\h_{\mu} = (\h^{\mu}, \circ)$
 this von Neumann algebra, and we  call
$\h_{\mu}$ the \emph{Poisson boundary} of $\mu$.

Our goal in this section is to prove quantum versions of several
important results which are well-known in the classical setting.
In order to prove our results in a general form for the Markov operators
corresponding to states on the universal $C^*$-algebra $\CU$, rather than just the ones in $\MG$,
we need to work with the universal von Neumann algebra $\CU^{**}$.
But there are some technical difficulties that arise in the non-Kac setting
if one wants to lift all quantum group properties to the universal von Neumann algebra (cf. \cite{Kus}).
So in the following we make sure that the properties we need for our purpose
are all valid at the universal von Neumann algebraic level.

Since $\LO$ is a norm closed two-sided ideal in $M_u(\G)$,
we obtain  a natural $M_u(\G)$-bimodule structure 
on $\LO$, and its adjoint defines 
an $M_u(\G)$-bimodule structure on $L_{\infty}(\G)$ such that 
\[
\langle f, \mu\star x\rangle = \langle f\star \mu, x\rangle \ \
~\mbox{and} ~ \ \
\langle f, x \star \mu\rangle = \langle \mu \star f, x\rangle
\]
for all $f\in \LO$, $\mu\in M_u(\G)$ and $x \in L_{\infty}(\G)$.
On the other hand, there is a natural  $M_u(\G)$-bimodule structure
on the universal enveloping von Neumann algebra $\CU^{**}$ given by
\[
\langle \omega, \mu\star_{u} x_{u}\rangle = \langle \omega\star_{u} \mu, x_{u}\rangle \ \
~\mbox{and} ~ \ \
\langle \omega, x_{u}\star_{u} \mu\rangle = \langle \mu \star_{u} \omega, x_{u}\rangle
\]
for all $\omega, \mu\in M_u(\G)$ and $x_{u} \in \CU^{**}$.
Let us denote by $\tilde \pi = (\pi^*\mid_{\LO})^{*} : \CU^{**}\to L_{\infty}(\G)$
  the normal  surjective $^*$-homomorphism extension of $\pi$ to $\CU^{**}$.
We obtain the following  interesting connection 
\begin{equation}
\label {F.full}
\tilde \pi(\mu\star_{u} x_{u}) = \mu \star \tilde \pi(x_{u})\ \
~\mbox{and }\ \
~ \tilde \pi(x_{u}\star_{u} \mu)= \tilde \pi(x_{u})\star \mu
\end{equation}
between the two module structures.
Indeed,  for any  $x_{u}\in \CU^{**}$,  we deduce from  
(\ref{F.2mul})  that
\begin{eqnarray*}
\langle f,  \mu\star \tilde \pi(x_{u}))\rangle =
\langle f \star \mu, \tilde \pi(x_{u}) \rangle
= \langle \pi^{*}(f) \star_{u} \mu, x_{u} \rangle
=  \langle \pi^{*}(f),  \mu \star_{u} x_{u} \rangle
= \langle f,  \tilde \pi(\mu \star_{u} x_{u}) \rangle
\end{eqnarray*}
for all $f\in\LO$ and $\mu\in M_u(\G)$.

The following result extends  \cite[Proposition 6.2]{Kus} to  the von Neumann algebraic level.
Here we denote by ${\mathcal V}$  the \emph{universal left regular corepresentation} 
of $C_{0}(\G)$ considered in \cite[Proposition 5.1]{Kus}.

\begin{proposition}\label {P.1.1}
For any $x_{u}\in \CU^{**}$ and $f \in \LO$ (or $f \in M(\G)$), we have 
\[
\pi^{*}(f)\star_{u} x_{u} = (\iota\otimes f){\mathcal V}^{*}(1 \otimes \tilde \pi(x_{u}))
{\mathcal V}.
\]
\end{proposition}
\begin{proof}
Given $x_{u}\in \CU^{**}$, there exists a net of elements $x_{i}\in \CU$ such that
$\|x_{i}\|\le \|x_{u}\|$ and $x_{i}\to x_{u}$ in the weak* topology.
It is known from 
\cite [Proposition 6.2]{Kus} that for each $x_{i}\in \CU$, we have 
\begin{equation}
(\iota \otimes \pi)\Gamma_{u}(x_{i}) = {\mathcal V}^{*}
(1 \otimes \pi(x_{i})){\mathcal V}.
\end{equation}
Therefore,  we get  
\begin{eqnarray*}
\langle \mu, \pi^{*}(f)\star_{u} x_{u}\rangle
 &=& \langle \mu \star_{u}\pi^{*}(f), x_{u}\rangle
 = \lim \langle \mu \star_{u}\pi^{*}(f), x_{i}\rangle\\
 &=& \lim \langle \mu \otimes \pi^{*}(f), \Gamma_{u}(x_{i})\rangle
= \lim \langle \mu \otimes f, {\mathcal V}^{*}(1 \otimes \pi(x_{i}))
{\mathcal V}\rangle\\
&=& \langle  \mu,  (\iota\otimes f){\mathcal V}^{*}(1 \otimes \tilde \pi(x_{u}))
{\mathcal V}\rangle
\end{eqnarray*}
for all   $\mu \in M_{u}(\G)$ and $f \in \LO$ (or $f \in M(\G)$).
\end{proof}

Using Proposition \ref{P.1.1}, we can prove the following result;
the idea of the proof is similar to the proof of \cite[Theorem 2.4]{Run2}.
\begin{proposition}\label{P.1.2}
For any $f\in \LO$ and $x_{u}\in \CU^{**}$, both $\pi^{*}(f) \star_{u}x_{u}$
and  $x_{u}\star_{u}\pi^{*}(f)$ are in $M(\CU)$.
\end{proposition}
\begin{proof}
For $f\in \LO$, we can write $f = y'\cdot f' $  for some $y'\in \mathcal{K}(L_{2}(\G))$
and $f' \in \LO$,
where $\mathcal{K}(\LT)$ denotes the $C^*$-algebra of all compact operators
on the Hilbert space $\LT$, and $\cdot$ is the canonical action of $\mathcal{K}(\LT)$
on its dual.
Since ${\mathcal V} \in M(\CU\otimes \mathcal{K}(L_{2}(\G)))$ (see \cite[Proposition 5.1]{Kus}), we have
\[
(\pi^{*}(f)\star _{u}x_{u})a 
= \langle \iota \otimes f', {\mathcal V}^{*} 
(1 \otimes \tilde \pi(x_{u})){\mathcal V}(a \otimes y')\rangle
\in \CU
\]
for all  $a\in \CU$. Here we used the fact that 
${\mathcal V}^{*} (1 \otimes \tilde \pi(x_{u})){\mathcal V}(a \otimes y')
\in \CU\otimes \mathcal{K}(L_{2}(\G))$.
This shows that $\pi^{*}(f) \star_{u} x_u\in M(\CU)$.
Similarly, we can prove that
$x_{u} \star_{u} \pi^{*}(f) \in M(\CU)$ by considering the
\emph{universal right regular corepresentation}
of $\CZ$.
\end{proof}

In the classical setting, when considering
the Poisson boundaries and harmonic functions on a locally compact group $G$,
in order to rule out trivialities,
one usually works with probability measures whose support generates $G$ as a closed
semigroup or group.
Therefore it is natural to seek for a quantum version of such a property and
restrict ourselves to those {\it quantum probability measures} possessing that property.

A state $\mu\in \PU_u(\G)$ is called  \emph{non-degenerate} on $\CU$
if for every non-zero element $x_u\in \CU^+$,
 there exists $n\in\field{N}$ such that
$\langle x_u , \mu^n\rangle \neq 0$ (see also \cite[Terminology 5.4]{VV3}).
Non-degeneracy can be  defined similarly for states  $\mu \in M(\G)$  on $\CZ$.
Note that every faithful state is non-degenerate, but there are examples of non-faithful
non-degenerate states.
If $\mu\in \PU_{u}(\G)$, then there exists a unique strictly continuous
state extension of $\mu$ to a state on $M(\CU)$, which we still denote by $\mu$.

\begin{lemma}\label{nondeg}
 Let $\mu\in \PU_u(\G)$ be non-degenerate. 
 Then for every non-zero $x_u\in M(\CU)^+$, there
exists $n\in\field{N}$ such that
$\langle x_u , \mu^n\rangle \neq 0$.
\end{lemma}
\begin{proof}
 Let $x_u\in M(\CU)^+$ be non-zero, and let $a_u\in\CU^+$ be such that 
 $\|a_u\| = 1$ and $a_{u}^{\frac 12}x_{u}^{\frac 12}\neq 0$. Then
we have $\CU^+\ni x_u^{\frac12}a_u x_u^{\frac12}\leq x_u$. Now since $\mu$ is non-degenerate,
there exists $n\in\field{N}$ such that
$0<\langle x_u^{\frac12}a_u x_u^{\frac12} , \mu^n\rangle\leq\langle x_u , \mu^n\rangle$.
\end{proof}

\begin{lemma}\label{P.faithful}
Let $\om\in \PU_u(\G)$,
and let $\psi$ and $\fee$ be, respectively, the right and the left Haar weights of $\G$.
Then $\Phi_{\om}$ is $\psi$-invariant and thus faithful on $L_{\infty}(\G)$;
similarly, the map $x\mapsto x\star \om$ is $\fee$-invariant, and hence faithful on $\LL$.
\end{lemma}
\begin{proof}
Since $\psi$ is the right Haar weight of $\G$, we have
\[
\psi (\Phi_{\om}(x))1=  (\psi\otimes \id)\Gamma (\Phi_{\om}(x)) 
= (\psi\otimes \Phi_{\om}) \Gamma(x)  =
\Phi_{\om}\left((\psi\otimes \id) \Gamma(x)\right)
= \psi(x)1
\]
for all $x \in \LL^{+}$.
This implies $\psi\circ \Phi_{\om}= \psi$ on $\LL^{+}$
and thus $\Phi_{\om}$ is faithful on $\LL$.
The result for the map $x\mapsto x\star \om$ follows similarly.
\end{proof}

The following lemma is essential for our results concerning non-degenerate states.
\begin{lemma}\label{lem2'}
Let $\G$ be a locally compact quantum group and let 
$\mu\in \PU_u(\G)$ be non-degenerate.
Let $x\in \LL$ be a self-adjoint element 
which attains its norm on $\LO_1^+$. If $x\in\mathcal{H}_\mu$
then $x\in\C1$.
\end{lemma}
\begin{proof}
Suppose that $\|x\| = 1$ and $f\in\LO^+$ is a state such that $\langle f , x \rangle = 1$.
Now assume towards a contradiction that $x\neq 1$.  Then $1-x$ is a non-zero positive element in
$\LL\cap \h_{\mu}$ and so there exists a non-zero positive element $x_u\in\CU^{**}$
such that $\tilde\pi(x_u) = 1-x$. Then by Proposition \ref{P.1.2}, we
have $x_u\star_u\pi^*(f)\in M(\CU)$. Moreover, by (\ref{F.full})
we have $\tilde\pi(x_u\star_u\pi^*(f)) = (1-x)\star f$.
It follows from Lemma \ref{P.faithful} (considering $\omega = f$ here) that $(1-x)\star f$ is non-zero,
which implies that $x_u\star_u\pi^*(f)\in M(\CU)$ is a non-zero
positive element. Since $\mu$ is non-degenerate, by Lemma \ref{nondeg},
there exists $n\in\field{N}$ such that
\[
\langle 1-x , f\star \mu^n\rangle =   \langle x_u\star_u\pi^*(f) , \mu^n\rangle \neq 0.
\]
On the other hand,  since $x\in\mathcal{H}_\mu$
we have $\Phi_{\mu^n}(x) = x$.
It follows that 
\[
 \langle 1 , f \star \mu^n\rangle = 1 = \langle x , f\rangle = \langle \Phi_{\mu^n}(x) ,
f\rangle = \langle x , f\star \mu^n\rangle.
\]
This implies that $\langle 1-x , f\star \mu^n\rangle = 0$, which is a contradiction.
Hence, we must have $x=1$.
\end{proof}

If  $\mu$ is a non-degenerate probability measure on a locally compact group $G$, it is
well-known that the space of all $\mu$-harmonic functions is never a subalgebra of $L_\infty(G)$, unless trivial.
Using the previous lemma, we can prove a quantum version of this result.
\begin{theorem}
 Let $\G$ be a locally compact quantum group and let
  $\mu\in \PU_u(\G)$ be  non-degenerate.
Then the following are equivalent:
\begin{itemize}
 \item [(i)] $\h_\mu$ is a subalgebra of $\LL$;
\item [(ii)] $\h_\mu = \C1$.
\end{itemize}
\end{theorem}
\begin{proof}
 We just need to prove (i) $\Rightarrow$ (ii). Since $\h_\mu$ is a weak$^*$ closed
operator system, (i) implies that $\h_\mu$ is a von Neumann subalgebra of $\LL$,
and is therefore generated by its projections. Now let $0\neq p\in\h_\mu$ be a projection
and $\xi\in\LT$ a unit vector such that $p\xi = \xi$.
Then we have $\|p\| = 1 = \langle p\xi , \xi\rangle$,
which shows that $p$ attains its norm on $\LO_1^+$.
Hence, $p=1$ by Lemma \ref{lem2'}. This shows that every projection of $\h_\mu$
is trivial and hence we have $\h_\mu = \C1$.
\end{proof}

It is also well-known that if $\mu$ is a non-degenerate measure on
a locally compact group $G$, then every continuous $\mu$-harmonic function on $G$
that vanishes at infinity is constant. We prove two non-commutative versions of this result.
\begin{theorem}
Let $\G$ be a locally compact quantum group and let $\mu\in \PU_u(\G)$ 
be non-degenerate. Then we have
$\mathcal{H}_\mu \cap \mathcal{K}(\LT) \subseteq \C 1$.
\end{theorem}
\begin{proof}
 It follows from the duality between $\mathcal{K}(\LT)$ and $\T(\LT)$, and the fact that
$\T(\LT)\mid_\LL = \LO$, that for every self-adjoint element $x\in\LL \cap \mathcal{K}(\LT)$,
either $x$ or $-x$ attains its norm on $\LO_1^+$. Hence, by Lemma \ref{lem2'}
we have $x\in\C 1$, and since $\mathcal{H}_\mu \cap \mathcal{K}(\LT)$ is generated by
its selfadjoint elements, the theorem follows.
\end{proof}

\begin{theorem}\label{CZ-tri}
Let $\G$ be a locally compact quantum group and 
let $\mu\in \PU_u(\G)$ be  non-degenerate. Then we have
$\mathcal{H}_\mu \cap \CZ \subseteq \C 1$.
\end{theorem}
\begin{proof}
Suppose that $x\in\mathcal{H}_\mu \cap \CZ$ is a 
self-adjoint element and that $\|x\| = 1$.
Then we can find (by substituting $x$ with $-x$, if necessary)
a state $\phi\in\MG = \CZ^*$ such that $\langle x , \phi \rangle = 1$.
Now, a similar argument to the proof of Lemma \ref{lem2'} shows that $x\in\C 1$.
Since $\h_\mu \cap \CZ$ is generated by its self-adjoint elements,
the theorem follows.
\end{proof}

As a consequence of Theorem \ref {CZ-tri}, we obtain the following result.

\begin{corollary}\label{0011}
Let $\G$ be a non-compact locally compact quantum group and
let  $\mu\in \PU_u(\G)$ be non-degenerate.
Then the Ces\`{a}ro sums $\displaystyle\{\frac{1}{n}(\mu + \mu^2 + \cdots + \mu^n)\}$
converge to $0$ in the weak* topology.
\end{corollary}
\begin{proof}
Let $\om\in M_u(\G)$ be an arbitrary weak$^*$ cluster point of
the Ces\`{a}ro sums $\displaystyle\{\frac{1}{n}(\mu + \mu^2 + \cdots + \mu^n)\}$
in $M_u(\G)$.
Then we get $\mu\star_u\om = \om$ and thus for any $x\in\CZ$, we have
\[
\Phi_\mu(\Phi_\om(x)) = \Phi_{\mu\star_u\om}(x) = \Phi_\om(x),
\]
which implies that $\Phi_\om(x)\in\h_\mu\cap\CZ$, and hence, by Theorem \ref {CZ-tri}, we have $\Phi_\om(x)\in\C 1$.
Since $\G$ is non-compact, this yields that $\Phi_\om(x) = 0$ for all $x\in\CZ$, and therefore it
follows from normality of the map $\Phi_\om$ that $\om=0$.
This shows that zero is the only weak$^*$ cluster point of
$\displaystyle\{\frac{1}{n}(\mu + \mu^2 + \cdots + \mu^n)\}$.
Since the unit ball of $M_u(\G)$ is weak$^*$ compact,
we get weak*$-\displaystyle\lim_n \frac{1}{n}(\mu + \mu^2 + \cdots + \mu^n) = 0$.
\end{proof}

The above Corollary \ref{0011} holds also for a non-degenerate state $\mu\in\MG$
(non-degenerate on $\CZ$), and as a consequence, we conclude the following.

\begin{corollary}\label{cpt-idem}
A locally compact quantum group $\G$ is compact if and only if there exists a
non-degenerate (and thus faithful) idempotent state $\mu\in M(\G)$.
\end{corollary}

Applying Corollary \ref {cpt-idem}, we obtain the following interesting result of
Fima \cite[Theorem 8]{fima}.

\begin{corollary}
Let $\G$ be a locally compact quantum group such that $\LL$ is a finite factor.
Then $\G$ is a compact Kac algebra.
\end{corollary}
\begin{proof}
Let $\tau\in\LO$ be the unique faithful trace. Then the uniqueness of $\tau$ implies that $\tau^2 = \tau$,
and hence, it follows from Corollary \ref{cpt-idem} that $\G$ is compact.
Moreover, it follows from \cite[Lemma 2.1]{Woro} that the trace $\tau$
is the Haar state of $\G$, and so $\G$ is a Kac algebra.
\end{proof}


\section{Amenability of Quantum Groups}

Our goal in this section is to prove a theorem establishing the equivalence
between amenability of a locally compact quantum group $\G$
and the absence of non-trivial harmonic operators on $\G$ (see Theorem \ref{KV}).
This answers the quantum group version of a conjecture formulated
in the group case by Furstenberg \cite{Furst}, which
in the classical setting was answered independently by
Kaimanovich--Vershik \cite{1} and Rosenblatt \cite{Ros}.

Let us first recall that a locally compact quantum
group $\G$ is \emph{amenable} if there exists a
\emph{left  invariant mean} on $\LL$, i.e., a state 
$F:\LL\to \C$
such that $(\id \otimes F)\Gamma (x) = F(x)1$.
Then a standard argument shows that  we can find a net  
of normal states $\{\omega_{\alpha}\}$
in $\loneqg$ such that 
\begin{equation}
\label {F.4equ}
\|f  \star\omega_{\alpha} - f(1)\omega_{\alpha}\|\to 0
\end{equation}
for all $f\in \loneqg.$
The following argument is inspired by  \cite [Theorem 4.3]{1}. An analogous result for the case of discrete Kac algebras was proved in \cite[Lemma 7.1]{Vaes}.
\begin{theorem} \label{P.4.2}
Let $\G$ be a locally compact quantum group such that $\loneqg$ is separable.
Then  the following are equivalent:
\begin{enumerate}
 \item  $\G$ is amenable;
 \item 
there exists a state $\mu\in\LO$ such that $\|f \star \mu^n - f(1)\mu^n \|
\,{\rightarrow}\, 0$ for
every $f\in\LO$, where $\mu^{n}= \mu \star \cdots \star \mu$ is 
the n-fold convolution of $\mu.$
\end{enumerate}
\end{theorem}
\begin{proof}
We only need to prove (1) $\Rightarrow $ (2).
  Let $\{f_i\}_{i\in\field N}$ be a dense subset of the unit ball of $\LO$, and 
  let $\{n_k\}$ be an increasing sequence of positive integers  such that
 $(\sum_{i=1} ^{k}\frac{1}{2^i})^{n_k} < \frac{1}{2^k}$.
Since $\G$ is amenable,  we can apply (\ref{F.4equ}) to  choose inductively 
a sequence  of states $\{\om_l\}_{l\in\mathbb N}$ 
in $\LO$ such that
\[
\| \om_{k_1}\star...\star\om_{k_r}\star \om_l - \om_l \| < \frac{1}{2^l}
\]
for all  $1\le k_i < l$ with $i= 1, \dots, r \le n_{l}$, and such that 
\[
\| f_s \star \om_{k_1}\star...\star\om_{k_r}\star \om_l - f_s(1)\om_l \|
 < \frac{1}{2^l}
\]
for all   $1\le s, k_i < l$ with $i= 1, \dots,  r\leq n_l$.
Define the normal state $\mu =  \sum_{l=1}^{\infty} \frac{1}{2^l} \om_l \in \loneqg$. 
Now given any  $f$ in the unit ball of $\LO$ and $\epsilon > 0$, we can choose $j\in\field N$ such that
$\| f - f_j \| < \epsilon$ and  $ \frac{1}{2^j} < \epsilon$. 
To simplify our notation, we   fix  $p= n_j$ and write  $t_i = \frac{1}{2^i}$. 
 Then we get 
\begin{eqnarray*}
\| f \star \mu^{p} - f(1)\mu^{p} \| &\leq& \|  f\star \mu^{p} - f_j \star \mu^{p} \|
 + \|f_j \star \mu^{p} - f_j(1)\mu^{p} \| + \| f_j(1)\mu^{p} - f(1)\mu^{p} \|\\
& <& 2\epsilon + \| f_j \star \mu^{p} - f_j(1)\mu^{p} \|.
\end{eqnarray*}
Now we split the term $\| f_j \star \mu^{p} - f_j(1)\mu^{p} \|$ as follows:
\begin{eqnarray*}
&&\| f_j \star \mu^{p}  - f_j(1)\mu^{p} \|= \| \sum_{max \ k_i\leq j} t_{k_1}...t_{k_p}
 f_j \star \om_{k_1}\star...\star\om_{k_p}
 + \sum_{max \ k_i > j} t_{k_1}...t_{k_p}  f_j \star \om_{k_1}\star...\star\om_{k_p}  \\
& - &\sum_{max \ k_i\leq j}  f_j(1) t_{k_1}...t_{k_p}\om_{k_1}\star...\star\om_{k_p}
 - \sum_{max \ k_i > j} f_j(1)  t_{k_1}...t_{k_p}\om_{k_1}\star...\star\om_{k_p} \| \\
&\leq& \sum_{max \ k_i\leq j} 2\|f_j\|  t_{k_1}...t_{k_p} +
 \sum_{max \ k_i > j} t_{k_1}...t_{k_p}  \| f_j \star \om_{k_1}\star...\star\om_{k_p} 
  -  f_j(1)\om_{k_1}\star...\star\om_{k_p} \|\\
&\leq& 2\epsilon +
 \sum_{max \ k_i > j} t_{k_1}...t_{k_p} \|  f_j \star \om_{k_1}\star...\star\om_{k_p}  
 - f_{j}(1)  \om_{k_1}\star...\star\om_{k_p}\|.
\end{eqnarray*}
Now consider one of the terms,
 $\om_{k_1}\star...\star\om_{k_p}$, in the last sum above
 and let $k_j$ be the smallest index such that $k_j > j$.
 Let $\mu_1 = \om_{k_1}\star...\star\om_{k_{j-1}}$ and $\mu_2 = \om_{k_{j+1}}\star...\star\om_{k_p}$. 
 Then we have
\begin{eqnarray*}
&&\|  f_j \star \om_{k_1}\star ...\star \om_{k_p}
- f_j(1)\om_{k_1}\star ...\star \om_{k_p} \| = 
\|   f_j \star \mu_1\star \om_{k_j}\star \mu_2   - f_j(1)\mu_1
\star \om_{k_j}\star \mu_2 \| 
\\ &\leq&\|  f_j \star  \mu_1 \star \om_{k_j}- f_j(1)  \mu_1 \star \om_{k_j}\| 
\leq \|  f_j \star  \mu_1 \star \om_{k_j} - f_j(1)\om_{k_j}\| \ + \ 
\| f_j(1) \mu_1 \star \om_{k_j} - f_j(1)\om_{k_j}\|
 < 2\epsilon,
\end{eqnarray*}
where the last inequality follows from the construction  of $\{\om_l\}$. This implies
\begin{eqnarray*}
\sum_{max \ k_i > j} t_{k_1}...t_{k_p}\|  f_j \star \om_{k_1}\star...\star\om_{k_p}
 - f_j(1)\om_{k_1}\star...\star\om_{k_p} \| < 2\epsilon .
\end{eqnarray*}
Hence we have $\|  f \star \mu^p - f(1)\mu^p \| < 6\epsilon$.
Since   $\|\mu\| = 1$,  we have 
\[
\|  f \star \mu^{p+l}- f(1)\mu^{p+l} \| \le \|  f \star \mu^{p}  - f(1)\mu^{p} \|
< 6\epsilon
\]
for all  $l \in \mathbb N$.
This implies that  $\| f \star \mu^n - f(1)\mu^n \| {\rightarrow} 0$.
\end{proof}


%

\begin{theorem}\label{KV}
Let $\G$ be a locally compact quantum group such that $L_{1}(\G)$ is separable.
Then the following are equivalent:
\begin{enumerate}
 \item  $\G$ is amenable;
\item there exists a state $\mu\in\MG$ such that $\h_{\mu}= \C1$.
\end{enumerate}
\end{theorem}
\begin{proof}
Recall that we denote by $\Phi_\mu$ the Markov operator $x\mapsto\mu\star x$ ($x\in\LL$).
Let us first assume that   $\G$ is amenable.  It follows from  Theorem \ref {P.4.2} that 
there exists $\mu\in\LO$ such that  
$\|f\star \mu^n  - f(1)\mu^n \| \to 0$
for every $f\in\LO$. 
Given any  $x\in {\mathcal H}_\mu$ and  $n\in \mathbb N$, we  have 
$\Phi_{\mu^{n}}(x) = \Phi_{\mu}^{n}(x)= x$.
It follows that for every $f\in \LO$, we have
\begin{eqnarray*}
\langle f , x - \mu^n(x) 1\rangle = \langle f , \Phi_{\mu^n }(x) - \mu^n(x) 1\rangle
= \langle  f \star \mu^n  - f(1)\mu^n , x\rangle\, {\rightarrow}\, 0.
\end{eqnarray*}
This  implies that $\mu^n(x) 1\to x$ in the weak* topology,
and thus we get  $x\in \C 1$.  This shows that $\h_{\mu}= \C1$.

On the other hand, let us suppose that we have a state $\mu\in M(\G)$ such that 
 $\mathcal{H}_\mu = \C 1$.
We choose a normal state $f\in\LO$. 
Then for each   $n\in \mathbb N$, we get a normal state
$\displaystyle\mu_n = \frac{1}{n}\sum_{k=1}^n  \mu^{k} \star f \in \loneqg$.
Let $F = \lim_{\mathcal U}\mu_n \in\LL^*$ be the weak* limit of $\{\mu_{n}\}$
with respect to a free ultrafilter  $\mathcal{U}$ on $\mathbb N$. 
Then $F$ is a state on $\LL$.
We claim that   $(\id \ot F)\Gamma(x)\in \h_{\mu}= \C1$ for all $x \in \LL$.
To see this, we notice that  the  Markov operator $\Phi_\mu$ satisfies 
\begin{eqnarray*}
\langle \Phi_\mu ((\id \ot F)\Gamma(x)), g \rangle   &=&
\lim _{\mathcal U}\frac{1}{n}\sum_{k=1}^n 
\langle \big (   \id \otimes \mu^{k}\star f  \big )  \Gamma(x), g \star \mu \rangle \\
&=&
\lim _{\mathcal U}\frac{1}{n}\sum_{k=1}^n
\langle (g\star \mu) \star (\mu^{k}\star f), 
x \rangle =  \langle (\id \ot F)(\Gamma(x)), g\rangle
\end{eqnarray*}
for all $g\in \loneqg$.
This shows that   $(\id \ot F)\Gamma(x)$ is an element in $ \mathcal{H}_{\mu}=\C 1$.

We define $F'\in\LL^*$ such that
$F'(x)1 = (\id \ot F)\Gamma(x)$.
Applying  $\mu$ to both sides of the latter,  we  obtain
\[
F'(x)= F'(x)\mu(1)= \mu (F'(x)1)= \mu((\id\otimes F)\Gamma(x))= F(x)
\]
for all $x \in \LL$.
Therefore $(\id \ot F')\Gamma(x)= F'(x)1 $ and thus $\G$ is amenable.
\end{proof}

\begin{remark} 
Using the same proof as that given in  Theorem  \ref {KV}, we can show 
that  a locally compact quantum group $\G$ is amenable if and only if
there exists a state $\omega\in \PU_u(\G)$ 
such that $\h_{\omega}= \C1$, if and only if there exists a  normal state $f\in\LO$
such that $\h_{f}= \C1$.
%
\end{remark}


\section{The Compact Quantum Group Case}

In this section, we consider compact quantum groups $\G$. Our goal is
to prove (Theorem \ref{cdthem}) a compact quantum group analogue of the
Choquet--Deny theorem.

Since $\G$ is compact, its reduced quantum group $C^*$-algebra $C_{0}(\G)$
and its  universal quantum group $C^*$-algebra  $C_{u}(\G)$
are unital Hopf $C^*$-algebras with the comultiplication
$\Gamma: C_{0}(\G)\to C_{0}(\G)\otimes C_{0}(\G)$
and the universal comultiplication $\Gamma_{u}: \CU\to \CU\otimes \CU$, respectively.
Also, in this case the $C^*$-algebra $C_0(\G)$ is
equal to the multiplier algebra $C(\G) = M(C_0(\G))$.
If $\phi$ is an idempotent state in $\PU_u(\G)$, i.e.,
 $\phi\star_u \phi= \phi$, it was shown in \cite [Theorem 4.1]{FS} that
$\tilde \h_{\phi}=\{x_{u}\in \CU: 
\tilde \Phi_{\phi}(x_{u})=x_{u}\}$ is a $C^*$-subalgebra of $\CU$, where
\begin{equation}
\label {F.idem}
\tilde \Phi_{\phi}(x_{u}) = (\id \otimes \phi)\Gamma_{u}(x_{u}) = \phi\star_u x_{u}.
\end{equation}
Using this fact, we can prove that for the corresponding Markov operator 
$\Phi_{\phi}= (\mathfrak{m}^{r}_{\phi})^{*}$ on $\LL$, 
the Poisson boundary $\h_{\phi}$ is a von Neumann subalgebra of $\LL.$
Let us first establish a lemma (see also \cite[Theorem 2.4]{SS}).

\begin{lemma}\label{l4.3}
Let $\G$ be a compact quantum group and let $\phi\in \PU_u(\G)$ be 
an idempotent state. Then we have
\begin{equation}\label{FS}
\Phi_{\phi}\big(\Phi_{\phi}(x) \Phi_{\phi}(y)\big)
= \Phi_{\phi}(x) \Phi_{\phi}(y)
\end{equation}
for all $x,y\in \LL$.
Moreover, the Poisson boundary $\h_{\phi}$ is a
von Neumann subalgebra of $\LL$.
\end{lemma}
\begin{proof} We first note that as an immediate consequence of (\ref {F.full}) and
(\ref {F.idem}), we get
\begin{equation}
\label {F.phiu}
\Phi_{\phi}(\pi(x_{u}))= \pi (\tilde\Phi_{\phi}(x_{u}))
\end{equation}
for all  $x_{u}\in \CU$.
Now for any $x, y\in C(\G)$, we can find $x_{u}, y_{u}\in \CU$ such that
$x = \pi(x_{u})$ and $y=\pi(y_{u})$, and thus we obtain
\begin{eqnarray*}
\Phi_{\phi}\big(\Phi_{\phi}(x) \Phi_{\phi}(y)\big) &=&
\Phi_{\phi}\big(\Phi_{\phi}(\pi(x_{u})) \Phi_{\phi}(\pi(y_{u}))\big)  =
\Phi_{\phi}\big(\pi(\tilde \Phi_{\phi}(x_u))\pi(\tilde \Phi_{\phi}(y_u))\big)\\
&=& 
\Phi_{\phi}\big(\pi(\tilde \Phi_{\phi}( x_u)\tilde \Phi_{\phi}(y_u))\big) =
\pi\big(\tilde \Phi_{\phi}(\tilde \Phi_{\phi}(x_u)\tilde \Phi_{\phi}(y_u))\big) 
\stackrel{(*)}{=}
\pi\big(\tilde \Phi_{\phi} (x_u)\tilde \Phi_{\phi}(y_u)\big) \\ &=&
\pi(\tilde \Phi_{\phi} (x_u))\pi(\tilde \Phi_{\phi} (y_u))=
\Phi_{\phi}(\pi(x_u))\Phi_{\phi}(\pi(y_u)) = \Phi_{\phi}(x)\Phi_{\phi}(y)
\end{eqnarray*}
where we used \cite[Theorem 4.1]{FS} in $(*)$.
It is known from the Kaplansky density theorem that the closed unit ball of $C(\G)$
is weak* dense in the closed unit ball of $\LL$.
Then for any contractive $x\in \h_{\phi}\subseteq \LL$, there exists a 
net of contractive elements  $x_{i}\in C(\G)$ such that $x_{i}\to x$ in the 
weak* topology.
Since $\Phi_{\phi}$ is weak* continuous, we get 
$\Phi_{\phi}(x_{i})\to  \Phi_{\phi}(x) = x$ in the weak* topology.
Similarly, for any $y\in \h_{\phi}$, we can find a net of elements 
$y_{j}\in C(\G)$ such that $\Phi_{\phi}(y_{j})\to y$ in the weak* topology.
Then we get the following iterated weak* limit 
\[
\Phi_{\phi}(x y)= \lim_{i}\lim_{j}\Phi_{\phi}(\Phi_{\phi}(x_{i})
\Phi_{\phi}(y_{j})) = 
\lim_{i}\lim_{j} \Phi_{\phi}(x_{i})
\Phi_{\phi}(y_{j}) = x y.
\]
This shows that the Choi--Effros product on $\h_{\phi}$ coincides with the product on
$\LL.$ Therefore, $\h_{\phi}$ is a von Neumann subalgebra of $\LL.$
\end{proof}

\begin{theorem}\label{sva2}
Let $\G$ be a compact quantum group and let $\mu$ be in $\PU_u(\G)$.
Then there exists an idempotent state $\phi\in \PU_u(\G)$ such that
$\h_{\mu}=\h_{\phi}$; in particular, $\h_\mu$ is a von Neumann subalgebra of $L_{\infty}(\G)$.
\end{theorem}
\begin{proof}
Consider the Ces\`{a}ro sums $\displaystyle\mu_n = \frac 1n(\mu + \dots + \mu^n)$, $n\in \mathbb{N}$,
and take the weak* limit $\phi = \lim_\mathcal{U}\mu_n$
with respect to a free ultrafilter $\mathcal U$ on $\mathbb N$.
Then  $\phi$ is an idempotent state in $ \PU_u(\G)$ 
such that $\phi\star_u \mu  = \phi = \mu \star_u \phi$.
Clearly, $\h_{\phi}\subseteq \h_{\mu}$ since for any
$x \in \h_{\phi}$, we have
\[
\Phi_{\mu}(x)= \Phi_{\mu}(\Phi_{\phi}(x))= \Phi_{\mu\star_u \phi}(x)
= \Phi_{\phi}(x)= x.
\]

To prove the converse inclusion, let us first suppose that 
$x\in C(\G)\cap\h_\mu$. Since $\G$ is compact, we have $C(\G) = \CZ$ and so there exists $x_u\in\CU$
such that $x = \pi(x_u)$. Then for any  $f\in\LO$, we have
\begin{eqnarray*}
\langle \Phi_\phi(x) , f\rangle &=&  \langle x , f \star \phi\rangle =
\langle x_u ,  \pi^*(f)\star_u \phi\rangle = \langle  x_u\star_u \pi^*(f) , \phi\rangle =
\lim_\mathcal{U}\langle  x_u\star_u \pi^*(f) , \mu_n\rangle \\
&=& \lim_\mathcal{U}\langle x_u , \pi^*(f)\star_u \mu_n\rangle =
\lim_\mathcal{U}\langle \Phi_{\mu_n}(x) , f\rangle 
= \langle x , f\rangle.
\end{eqnarray*}
This shows that $x\in\h_\phi$. Hence, $\h_\mu\cap C(\G)\subseteq\h_\phi$.

Now let  $x\in\h_\mu$. Then for any $f\in\LO$ we have  $ x\star f \in LUC(\G) 
= C(\G)$.  We also have   $x \star f \in \h_{\mu}$ since
$\Phi_{\mu}( x\star f ) = \Phi_\mu(x)\star f = x \star f$.
Therefore we have $x\star f\in C(\G)\cap\h_\mu\subseteq\h_\phi$
for all $f\in\LO$.
From this we conclude that
\[
\langle \Phi_\phi(x) , f\star g\rangle = \langle x , f\star g\star\phi\rangle =
\langle x\star f, g\star\phi\rangle = \langle \Phi_\phi(x\star f) , g\rangle
=  \langle x\star f , g\rangle =  \langle x , f\star g\rangle
\]
for all $f, g \in \LO$.
Since $\LO= \langle \LO\star\LO \rangle$,
we obtain $\Phi_\phi(x) = x$. Hence, $\h_\mu\subseteq\h_\phi$.
\end{proof}

Now, as a corollary to Theorems \ref{CZ-tri} and \ref{sva2}, we have the
following compact quantum group analogue of the Choquet--Deny theorem.
A special case of this result was proved by Franz and Skalski \cite{FS2}
where $\mu$ was assumed to be faithful.


\begin{theorem}\label {cdthem}
Let $\G$ be a compact quantum group and let $\mu\in \PU_u(\G)$
be non-degenerate.
Then we have $\mathcal{H}_\mu=\C 1$.
\end{theorem}

\section{Examples}

It is often highly non-trivial
to concretely identify Poisson boundaries associated to
a given locally compact quantum group.
The situation in the classical setting is of course much easier.
The structure of Poisson boundaries 
has been studied in detail for locally compact groups
in many interesting cases.

In this section we will establish a bridge between the classical and the quantum setting,
through a concrete formula (\ref{cla-qua-for}),
which then allows us to link the Poisson boundaries in these two settings.
 In particular, we apply our machinery to the
case of Woronowicz' twisted $SU_q(2)$ and show that the 
Poisson boundary associated to a specific state on this compact quantum
group can be identified with the Podle\'{s}  sphere.

Throughout this section, $\G$ denotes a co-amenable locally compact quantum group.
Let us recall that in this case we have
$M(\G) = M_u(\G)$.

\begin{proposition}\label{Pois-dense}
 Let $\G$ be a co-amenable locally compact quantum group.
 If  $\mu$ is a state in  $M(\G)=M_u(\G)$, then the closed unit ball of
  $\h_{\mu}\cap LUC(\G)$ is
 weak$^*$ dense in the closed unit ball of  $\h_\mu$.
\end{proposition}
\begin{proof}
Let $y\in\h_\mu$ with $\|y\| \le 1$ and let $\{f_\alpha\}\subseteq\LO$ be 
a contractive approximate identity.
  Then we have
$$\Phi_\mu\big((f_\alpha\ot\id)\Gamma(y)\big) 
= (f_\alpha\ot\id)(\id\ot\Phi_\mu)\Gamma(y) =
(f_\alpha\ot\id)\Gamma(\Phi_\mu(y)) = (f_\alpha\ot\id)\Gamma(y),$$
which implies that $y \star f_{\alpha}=(f_\alpha\ot\id)\Gamma(y)
\in\h_\mu\cap LUC(\G)$.
Since $\{f_\alpha\}$ is a contractive approximate identity, we 
have $\|y \star f_{\alpha}\|\le 1$ and 
$(f_\alpha\ot\id)\Gamma(y)\to  y$ in the weak* topology.
This  completes the proof.
\end{proof}


It  was shown in  Kalantar's thesis \cite  [Chapter 3]{mehrdad-thesis}  that 
if $\G$ is a co-amenable locally compact quantum group, the \emph{spectrum} 
\[
\tilde\G = sp(\CZ) = \{\phi: C_{0}(\G) \to \C \ | \ \phi ~\mbox{is a non-zero *-homomorphism}\}
\]
of $\CZ$  equipped with the convolution product 
and the weak$^*$ topology from $\MG$  is a locally compact group.
We let
\[
\wedge : x \in C_{0}(\G) \to \hat x \in  C_{0}(\G)^{{**}}
\]
be the canonical second dual inclusion and let 
\begin{equation}
\label {F.P}
P: x \in C_{0}(\G) \to \hat x\mid_{\tilde \G}\in C_{0}(\tilde \G)
\end{equation}
be the Gelfand transformation given by $P(x)(\phi) = \phi(x)$ for all $\phi \in \tilde \G$.

\begin{proposition}\label {F.Kalantar}
Let $\G$ be a co-amenable locally compact quantum group.
The map $P$ defined in (\ref {F.P})
is a *-homomorphism from $C_{0}(\G)$ onto $C_{0}(\tilde \G)$,
and thus has a unique strictly continuous unital *-homomorphism extension from 
the $C^*$-multiplier algebra $C(\G)= M(C_{0}(\G)) $ onto the 
$C^*$-multiplier algebra   $C(\tilde \G)= M(C_{0}(\tilde \G))$.
\end{proposition}
\begin{proof} It is easy to see that  $P$ is a *-homomorphism from
$\CZ$ into $C_{0}(\tilde \G)$.
Then the range space of $P$  is a $C^*$-subalgebra of  $C_{0}(\tilde \G)$ and it 
separates points in $\tilde \G$.
Therefore, by the generalized Stone-Weierstrass theorem, 
we have $P(\CZ) = C_0(\tilde\G)$.
Therefore,   $P$ has a unique  strictly continuous  unital *-homomorphism extension,
mapping the $C^*$-multiplier algebra $C(\G)= M(C_{0}(\G))$ onto the 
$C^*$-multiplier algebra 
 $C(\tilde \G)= M(C_{0}(\tilde \G))$ (cf. Lance \cite {Lance}).
\end{proof}

Since the comultiplication $\Gamma:  \CZ \to M(\CZ\otimes \CZ)$ 
and the *-homomorphisms
\[
\id \otimes P : \CZ \otimes \CZ \to \CZ\otimes C_{0}(\tilde \G)~ 
\mbox{and} ~ 
P \otimes P: \CZ \otimes \CZ \to C_{0}(\tilde \G)\otimes C_{0}(\tilde \G)
\]
 have  unique  strictly continuous
extensions   to their $C^*$-multiplier algebras, we can consider their 
compositions with the comultiplication
$\Gamma$ and obtain the following result.
\begin{proposition}\label {P.P} Let $\G$ be a co-amenable locally compact quantum group
and let $P: C(\G)\ra C(\tilde\G)$
 be the  strictly continuous  unital *-homomorphism defined above.
 \begin{itemize}
 \item [(1)] If we let  $\Gamma_{a}$ denote the comultiplication on    $C_0(\tilde \G)$,
 we have 
 \begin{equation}\label {F.comm}
 (P\otimes P)\circ \Gamma = \Gamma_{a}\circ P.
 \end{equation}
 \item [(2)] The induced map $(\id \otimes P)\circ \Gamma$ is an  injective *-homomorphisms from 
 $C(\G)$ into $M(C_{0}(\G) \otimes C_{0}(\tilde \G))$.
\end{itemize}
\end{proposition}
\begin{proof}
The first part follows from straightforward calculations. For the second part,
let $\varepsilon$ be the unital element in $M(\G)$.
So, we have $(\id \otimes \varepsilon)\Gamma(x) = x$ for all $x \in \CZ$,
 and thus for all $x \in C(\G)$.
 Moreover, $\varepsilon$ is a non-zero *-homomorphism,
and thus is an element, which is denoted by $e$, in $\tilde \G$. Since the
multiplication of the group $\tilde\G$ is induced from the multiplication of $M(\G)$, $e$ is just
the unital element of $\tilde\G$. Moreover, for any $x\in \CZ$, we have
\[
\varepsilon(x) = \hat x(\varepsilon) = P(x)(e) = e(P(x)).
\]
This implies that $\varepsilon = e \circ P$.
Now, if we are given  $x \in C(\G)$ such that  $(\id\ot P) \Gamma (x) = 0$, then we have
\[
 x = (\id\ot \varepsilon)\Gamma (x) =  (\id\ot e)(\id \otimes P)\Gamma (x) = 0.
\]
So $(\id \otimes P)\circ \Gamma $ is injective.
\end{proof}

Since $P$ is a  *-homomorphism from $\CZ$ onto $C_{0}(\tilde \G)$ satisfying 
(\ref {F.comm}), its adjoint map $P^*$  defines  a completely  isometric and 
Banach algebraic homomorphism 
\begin{equation}
\label {F.P*}
P^{*}:  M(\tilde\G)\ni\mu\mapsto \mu_{\G}=  \mu \circ P\in\MG.
\end{equation}
Therefore, we can identify $M(\tilde \G)$ with a Banach subalgebra of  $\MG$.

\begin{theorem}\label{qcident}
 Let $\G$ be a co-amenable locally compact quantum group and let $\mu$ be a 
 probability  measure in $M(\tilde\G)$. Then we have
\begin{equation}\label{cla-qua-for}
 \h_{\mu_\G} = 
 \overline{\{x\in C(\G) : (\id\ot P)\Gamma(x) \in \LL\vtp\h_{\mu}\}}^{\text{weak}^*}.
\end{equation}
\end{theorem}
\begin{proof}
We prove that
\[
 \h_{\mu_\G}\cap C(\G) = \{x\in C(\G) : (\id\ot P)\Gamma(x)
 \in \LL\vtp\h_{\mu_\G}\}.
 \]
from which the result follows by Proposition \ref {Pois-dense}.

  Given any $x \in C(\G)$, we have by (\ref{F.comm}) that
\[
 \Phi_\mu\circ P(x) = (\id\ot\mu)\Gamma_a(P(x)) = (\id\ot\mu)(P\ot P)(\Gamma(x))
= P\big((\id\ot\mu_\G)\Gamma(x)\big) = P\circ\Phi_{\mu_\G}(x).
\]
 Therefore, if $x\in\h_{\mu_\G}\cap C(\G)$, then we get 
 $P(x)\in\h_{\mu}\cap C(\tilde\G)$ and 
\[
 (\id\ot \Phi_\mu)(\id\ot P)\Gamma(x) 
 = (\id\ot P)(\id\ot\Phi_{\mu_\G} )\Gamma(x) =
(\id\ot P)\Gamma(\Phi_{\mu_\G}(x)) = (\id\ot P)\Gamma(x).
\]
Then, one can show that $(\id\ot P) \Gamma(x)\in\LL\vtp\h_{\mu}$ (cf. \cite{KNR2}).

On the other hand, assume that $x\in C(\G)$ is such that 
$(\id\ot P)\Gamma(x)\in\LL\vtp\h_{\mu}$.
Then we have
\[
(\id\ot P)\Gamma(\Phi_{\mu_\G}(x)) = (\id\ot P)(\id\ot\Phi_{\mu_\G} )\Gamma(x) =
  (\id\ot \Phi_\mu)(\id\ot P)\Gamma(x)  =  (\id\ot P)\Gamma(x).
\]
Hence, by Proposition \ref{P.P}, we obtain $\Phi_{\mu_\G}(x) = x$.
\end{proof}

In the following, we consider Woronowicz's $SU_q(2)$ quantum group
for $q\in (-1, 1)$ and $q\neq 0$ (cf. \cite{Wosuq}).
It is known that $SU_q(2)$ is a co-amenable compact quantum group 
with the quantum group $C^*$-algebra $C(SU_q(2))= C_{u}(SU_{q}(2))$ generated by 
two operators $u$ and $v$ such that   $U  = \left [ \begin{array}{cc}
  u & -qv^{*}  \\
  v & u^{*} 
\end{array} \right] $   is a unitary matrix in $M_{2}(C(SU_q(2)))$.

It was shown in \cite [Theorem 3.4.3]{mehrdad-thesis} that 
$\tilde \G$ is actually homeomorphic to  the unit circle group $\mathbb T$.
Indeed, if  $f \in \tilde \G$ is a non-zero *-homomorphism 
on $C(SU_q(2))$,  then 
\[
\left [ \begin{array}{cc}
 f(u) & f(-qv^{*})  \\
  f(v) & f(u^{*}) 
\end{array} \right] 
 = \left [\begin{array}{cc}
  f(u) & -q \overline {f(v)}  \\
  f(v) & \overline {f(u)}
\end{array} \right] \ 
\]
is a unitary matrix in $M_{2}(\C)$.  This implies that 
\[
|f(u)| ^{2} + |f(v)|^{2} = 1 ~\mbox{and} ~ |f(u)| ^{2} + q^{2}|f(v)|^{2} = 1.
\]
Since $0 < |q| < 1$, we must have $f(v)= 0$ and $|f(u)| = 1$.
Then we  get a map 
\[
\gamma: \widetilde{SU_q(2)}\ni f \mapsto \left [ \begin{array}{cc}
 f(u) & 0  \\
  0 & \overline {f(u) }
\end{array} \right]
\]
which gives a map
from $\widetilde {SU_q(2)}$ into the unit circle $\mathbb T$.
Since $C(SU_q(2))$ is the universal $C^*$-algebra 
generated by $u$ and $v$, 
 it is easy to see that 
$\gamma$ defines  a homeomorphism  from $\widetilde{SU_q(2)}$ onto 
$\mathbb T$.
Moreover, since $\widetilde{SU_q(2) }$ is a compact group, 
and $\Gamma(u) = u \ot u$
 (see \cite [Theorem 1.4]{Wosuq}),  $\gamma$ is a group
homeomorphism from $\widetilde{SU_q(2)}$ onto $\mathbb T$.
Therefore, we can identify the spectrum $\widetilde{SU_q(2)}$ with $\mathbb T$.

In view of the above discussion, the *-homomorphism $P$ defined in (\ref {F.P})
can be identified with a map 
 \[
P_{\gamma}: x \in  C(SU_q(2))\ra \hat x \circ \gamma^{-1} \in C(\mathbb{T})
\]
such that $P_{\gamma}(u)= \id_\mathbb{T}$ and $P_{\gamma}(v)  = 0$,
where $\id_\mathbb{T}:\mathbb{T}\ra \mathbb{T}$ is the identity function 
$z\mapsto z$ on $\mathbb{T}$.
Now  let
\begin{equation}\label{qsph}
 C(SU_q(2)\backslash\mathbb{T}) = 
 \{x\in C(SU_q(2)) : (\id\ot P_{\gamma})\circ \Gamma(x) = x\ot1\}.
\end{equation}
Then $C(SU_q(2)\backslash\mathbb{T})$ is a $C^*$-subalgebra of $C(SU_q(2))$, 
and one can show that $(C(SU_q(2)\backslash\mathbb{T}),
\Gamma\mid_{C(SU_q(2)\backslash\mathbb{T})})$
is one of the Podles' quantum spheres (see \cite{Podles} for the details).
We also call the von Neumann algebra generated
by $C(SU_q(2)\backslash\mathbb{T})$ in $L_\infty(SU_q(2))$ a
\emph{quantum sphere} and will denote it by  $L_\infty(SU_q(2)\backslash\mathbb{T})$.
In the next theorem, we show that the quantum sphere 
$L_\infty(SU_q(2)\backslash\mathbb{T})$
is a concrete realization of the Poisson boundary of Markov operators associated with
non-degenerate measures in $M({\mathbb T})$.

\begin{theorem}
 Let $\mu\in M({\mathbb T})$ be a non-degenerate  measure. Then we have 
 \begin{equation}
 \h_{\mu_{SU_q(2)}}= L_\infty(SU_q(2)\backslash\mathbb{T}).
 \end{equation}
\end{theorem}
\begin{proof}
It follows from (\ref{F.P*}) and (\ref {qsph}) that if 
$x \in C(SU_q(2)\backslash\mathbb{T})$,
then we have 
\[
\Phi_{\mu_{{SU_{q}(2)}}} (x) = (\id \otimes \mu) (\id \otimes P_{\gamma})
\Gamma(x) =  (\id \otimes \mu)(x \otimes 1)= x.
\]
Hence we see that $L_\infty(SU_q(2)\backslash\mathbb{T})\subseteq \h_{\mu_{SU_q(2)}}$.
For the converse inclusion, by Theorem \ref{qcident} it is enough
to show that any $x\in C(SU_q(2))$ with
$$(\id\ot P_{\gamma})\Gamma(x) \in \LL\overline {\ot} \h_\mu$$
lies in $C(SU_q(2)\backslash\mathbb{T})$.
 Since $\mu\in M(\mathbb{T})$ is non-degenerate, $\h_\mu = \C1$,
 hence $(\id\ot P_{\gamma})\Gamma(x) \in \LL\overline {\ot} \C 1$.
So, there exists $y\in\LL$ such that $(\id\ot P_\gamma)\Gamma(x) = y\ot 1$.
This implies that
\[
 x = (\id\ot\varepsilon)\Gamma(x) = (\id\ot e)(\id\ot P) \Gamma(x) 
 = (\id\ot e)(y\ot 1) = y,
\]
which yields that $x\in C(SU_q(2)\backslash\mathbb{T})$.
\end{proof}




\bibliographystyle{plain}

\end{document}